\documentclass{article}

\usepackage{amsmath}
\usepackage{amsfonts}
\usepackage{amssymb}
\usepackage{amsthm}
\usepackage{graphics}
\usepackage{enumerate}
\usepackage{color}
\usepackage{tikz-cd}[2010/10/13] 
\usepackage[bookmarks=false]{hyperref}
\usepackage{cleveref}
\usepackage{tensor}
\usepackage{booktabs}
\usepackage[normalem]{ulem}
\usetikzlibrary{arrows}

\newcommand{\D}{\nabla}
\newcommand{\ed}{\mathrm{d}}

\newcommand{\R}{\mathbb{R}}
\newcommand{\E}{\mathbb{E}}

\theoremstyle{definition}

\newtheorem*{remark*}{Remark}
\DeclareMathOperator{\Image}{Im}

\title{Optimal Projection Filters\thanks{An updated version of this preprint will appear in the journal Information Geometry, for the special issue ``Half a century of Information Geometry"}}
\author
{Damiano Brigo \\ Dept.\ of Mathematics \\ Imperial College London \\ {\small \tt{damiano.brigo@imperial.ac.uk}}}

\date{First version: Feb. 21, 2022. This version: May 3, 2022}

\begin{document}

\maketitle
\begin{abstract}

We present the two new notions of projection of a stochastic differential equation (SDE) onto a submanifold, as developed in Armstrong, Brigo e Rossi Ferrucci (2019, 2018) \cite{armstronglms,rossiferrucci}: the It\^o-vector and It\^o-jet projections. This allows one to systematically and optimally develop low dimensional approximations to high dimensional SDEs using differential geometric techniques.  Our new projections are based on optimality arguments and yield a well-defined ``optimal'' approximation to the original SDE in the mean-square sense. We also show that the earlier Stratonovich projection satisfies an optimality criterion that is more ad hoc and less natural than the criteria satisfied by the new projections. As an application, we consider approximating the solution of the non-linear filtering problem within a given manifold of densities, using either the Hellinger or $L^2$ direct metrics and related Information Geometry structures on the space of densities. The Stratonovich projection had yielded the projection filters studied in Brigo, Hanzon and Le Gland (1998, 1999) \cite{brigo1,brigo2}, while the new projections lead to the optimal projection filters. The optimal projection filters have been introduced in \cite{armstronglms}, where numerical examples for the Gaussian case are  given and where they are compared to more traditional nonlinear filters.   
\end{abstract}

\bigskip


\bigskip

{AMS classification codes: 62M20, 93E11, 60G35, 62B10, 58J65, 60H10, 65D18, 58A20}
\newpage

\thispagestyle{empty}

\tableofcontents

\newpage

\section{Introduction and history}

\subsection{Information geometry as the differential geometric approach to statistics}

Information Geometry is an informal term to describe the differential geometric approach to statistics, or more precisely to study the differential geometric properties of sets of probability distributions, on which a manifold structure is usually built, leading to so called statistical manifolds. The beginning of information geometry is usually attributed to the legendary mathematician and statistician, C. R. Rao, starting with his 1945 paper \cite{rao45}. The theory has as one of its originating points the interpretation of the Fisher matrix for a parametric family of distributions as a Riemannian metric on the given finite-dimensional statistical manifold, the dimension being usually related to the number of parameters. The Fisher information matrix is related naturally to the Hellinger distance on more general infinite-dimensional spaces of probability measures, a distance based on the $L^2$ structure of sets of square roots of probability densities.  Nand Lal Aggarwal (1974)\cite{aggarwal},  Shun\textquotesingle ichi Amari (1985) \cite{amari},  Ole Barndorff-Nielsen (1978) \cite{barndorff} and   Giovanni Pistone (Pistone and Sempi 1995 \cite{pistoneannals}) are key initial references, among others, that contributed to the development of information geometry. 

\subsection{Information geometry and filtering dynamics}

Our work concerns the application of information geometry to approximation of dynamics of probability distributions, in most cases  stemming from the stochastic filtering problem. 

To state it in basic terms, in stochastic filtering one observes a random signal perturbed by random noise. The unperturbed random signal cannot be observed but needs to be estimated. For example, the perturbed signal could be the radar reading of the position of a spacecraft, which would not provide the exact position of the spacecraft due to several disturbances (``noise'') in the radar observations. It would then be necessary to estimate the real position of the spacecraft from the noisy radar readings. This is a filtering problem. A filtering algorithm was used in the Apollo 11 mission (Cipra 1993 \cite{cipra}), the first human landing on the moon. Filtering has also applications in areas such as water level estimation and prediction, submarine navigation, econometrics, target tracking and many others. A good historical book on filtering with an eye to applications is Jazwinski (1970) \cite{jazwinski}, see also Maybeck (1982) \cite{maybeck}, while the mathematical aspects are considered fully in Liptser and Shiryayev (1978) \cite{liptser}. More recent monographs on filtering are Ahmed (1998) \cite{ahmedbook} and Bain and Crisan (2009) \cite{crisan}.

The general solution of the filtering problem at a given time is given by the probability density of the unperturbed state of the system at that time, conditional on the perturbed observations up to the given time. When the unobserved signal and the observed signal evolve in continuous time, the filter density follows a stochastic partial differential equation (SPDE). It has been shown that this probability density, the solution of the SPDE, does not evolve in a finite dimensional statistical manifold, except in very special cases. For example, if the dynamics of the unobserved system is linear, the observations are linear, the noises are Gaussian and the initial condition of the unperturbed signal is also Gaussian (or deterministic), then the filter is Gaussian and its density can be characterized by a finite dimensional set of parameters, namely the mean vector and variance-covariance matrix of the resulting Gaussian distribution. This leads to the celebrated Kalman filter. However, this does not happen usually, in the non-linear case, and the filter is infinite dimensional in general, as shown for the cubic sensor example by Hazewinkel, Marcus and Sussmann (1983) \cite{hazewinkel}.

\subsection{Classic projection filters: Stratonovich--Hellinger projection}\label{sec:clastratpr}

Enters information geometry. Can information geometry provide us with a method to approximate the infinite-dimensional filter with a finite-dimensional approximation that is close to the original filter? The idea to apply the Fisher Metric and Hellinger distance to this problem was first sketched in an article of Bernard Hanzon (1987) \cite{hanzon} while he was working at the Technical University of Delft. Hanzon suggested to project the SPDE equation in Stratonovich form for the evolution of the filter density onto a finite dimensional statistical manifold, using the Fisher metric/Hellinger distance. We call this ``Stratonovich projection'' and it consists in projecting the separate vector fields of the SPDE corresponding to the drift and diffusion part of the Stratonovich version. The projected equation would describe a finite dimensional density evolution, called projection filter, approximating the full filter evolution associated with the optimal filter. The paper was presented to a key conference in Lancaster whose proceedings were edited by  Christopher T. J. Dodson, in a volume with the almost prophetic title ``Geometrization of Statistical Theory''. The following year, on August 22, 1988, Hanzon presented the idea at a seminar in Tokyo University called ``The Projection Filter'' while visiting Shun\textquotesingle ichi Amari. A few years later, in 1991, Hanzon and a PhD student Ruud Hut also from Technical University of Delft, wrote the paper Hanzon and Hut \cite{hanzonhut} with new results on the projection filter on Gaussian densities, showing that for the Gaussian family the projection filter coincides with a heuristic-based family of finite dimensional filters, the assumed density filters, previously studied by Harold Kushner (1967) \cite{kushner}, see also \cite{maybeck}.    

The projection filter idea was formulated precisely, extended and made fully rigorous in subsequent works, during the PhD studies of Damiano Brigo with Bernard Hanzon at the Free University of Amsterdam and with Francois LeGland at IRISA/INRIA, in Rennes, France, in 1993--1996 \cite{brigoPHD}. In these studies it was shown, among other things, that exponential families played a very particular role in the projection filter, allowing for the correction step of the filtering algorithm to be exact, and also fully generalizing the equivalence to the assumed density filters. The filters were tested numerically on some examples. During his PhD, Brigo also authored other papers on small observation noise for the Gaussian projection filter \cite{brigoscl1,brigoscl2} and on approximations of the Fokker-Planck-Kolmogorov equation, as well as formulations of the filter in discrete time using the Kullback Leibler information, with application to volatility modeling in finance \cite{brigoIME}. The main results on the projection filters were published later in Brigo, Hanzon and LeGland (1998, 1999) \cite{brigo1,brigo2}.

One of the key issues, from the start, was making sure that the given approximated equation for the filter density would stay on the chosen statistical manifold. The Stratonovich projection  ensured this, but scholars had been studying the behaviour of stochastic differential equations on manifolds independently of the filtering application above. Among those, we refer to  David Elworthy (1988) \cite{elworthy},  Michel Emery (1989) \cite{emery}, and more recently  Elton Hsu (2002) \cite{hsu}. We also notice that Elworthy, Le Jan and Li (2010) would discuss geometric aspects of filtering theory in \cite{geometryoffiltering}, although their book does not deal with projection filters. 

\subsection{Classic projection filters: Stratonovich--direct $L^2$ projection}\label{sec:classtratl2}

The cold war had ended, the military applications of filtering were receiving less and less funding, and many stochastic analysts who had previously worked on filtering turned to other areas, and many turned to mathematical finance. Brigo turned to a career in mathematical finance, but returned to the filtering problem as a side project in 2011 after he moved from a managing director position in the financial industry to a full academic position as Gilbart Chair at the Department of Mathematics of King's College London, earlier in 2010. There, in 2011 he met a new colleague, John Armstrong, a differential geometry PhD from Oxford who had worked on almost K{\"{a}}hler geometry and who also had spent several years in the financial industry and was now turning to a full-time academic career. Brigo explained the filtering problem to Armstrong, who grasped immediately the essential ideas and the mathematics. Brigo had already written a preprint on his new idea of applying the direct $L^2$ structure without square roots to obtain a new type of projection filter, showing equivalence with Galerkin-based filters when using mixtures of distributions. Armstrong refined the idea and implemented the filter numerically, studying the cubic sensor problem. This led to a second wave of projection filters based on the direct $L^2$ metric as opposed to the Hellinger distance. It turned out that, as anticipated in the preprint, while the original Hellinger-based filters worked well with exponential families, being equivalent to assumed density filters, the direct $L^2$ filters worked best with mixture families, being equivalent to Galerkin-based filters. This research went on in 2011-2013 and was published in Armstrong and Brigo (2016) \cite{armstrongbrigomcss}. By 2012 Brigo had moved to Imperial College, so this had become a cooperation between the Mathematics  Departments at King's College London and Imperial College London. During the review of \cite{armstrongbrigomcss}, one of the reviewers asked in which sense, or according to which criterion, the projection filter was providing an optimal approximation of the true filter. Armstrong immediately grasped the problem, while Brigo realized that while he had taken the Stratonovich projection for granted, he could not really tell in which sense it was optimal for the approximation of the true filter SPDE as a whole.  

\subsection{Is the classic projection filter an optimal approximation?}

The essence of the problem of optimality of the approximation was based on the way the filtering equation was projected in the projection filter works published until then, mainly   \cite{brigo1,brigo2,armstrongbrigomcss}.  There are two stochastic calculi, Ito and Stratonovich. The two different calculi are suited to different applications, but from a probabilistic point of view the Ito calculus has a more clear interpretation of the stochastic equation coefficients in terms of local mean and local standard deviation, linked to the martingale property. Also, it is believed that even when one works with Stratonovich calculus, under the formalism one can argue that it is still the Ito calculus
that ``does all the work'' (Rogers and Williams (1987) \cite{rogerswilliams}, Chapter V.30, p. 184).
The problem with Ito calculus is that it violates the chain rule for change of variables. When changing variables, one has to use the famous Ito's formula, involving a second order term in the transformation. 

The true, infinite dimensional filter equation (taking the form of a stochastic partial differential equation, or SPDE) had always been written in Stratonovich form in the previous projection filter works, because in a Stratonovich stochastic equation the two parts describing the drift term and the diffusion coefficient term obey the chain rule under change of variables. This means that they can be interpreted as vector fields and be projected without problems on the tangent space of a submanifold, obtaining vector fields in the submanifolds that would form the approximating finite dimensional stochastic differential equation. 

Projecting directly the Ito equation does not work, because the change of variables includes second order terms that do not resemble the behaviour of vector fields. Projection becomes then impossible to perform directly in Ito form. One could re-write the Ito true filter stochastic equation in Stratonovich form, project it, obtain a finite dimensional approximated filter, and transform back this approximate filter equation from Stratonovich to Ito form. But in what sense is this approximation optimal? What criterion does it minimize? 
 
Brigo had never thought about this in depth because, in the back of his mind, he believed the projection of a vector field to always provide the best approximation of the original vector field. But a stochastic equation is given by two terms, the drift and the diffusion part, and if one puts the equation in Stratonovich form, the drift and the diffusion coefficients become described by two vector fields and as such can be projected. As the two vector fields are projected, each projected vector field will be  the best approximation of the original vector field, but what does this mean for the solution of the stochastic equation as a whole? The stochastic equation is not just the pair of vector fields. In fact, when the equation is in Ito form, the drift and the diffusion coefficients interact when changing variables or coordinates, involving second order terms in the transformation. The fact that Stratonovich is ``less good'' probabilistically means that putting together two optimal projections of the coefficients to form a single Stratonovich equation does not provide a solution that is optimal in a probabilistic sense, for example in mean square.

\subsection{Finding optimal projection filters}
Armstrong had previously noticed that an Ito equation behaved exactly as a geometric object he was familiar with, called a 2-jet. Brigo, while helping Armstrong in developing the 2-jet interpretation of stochastic differential equations, started looking at the two other known ways to model stochastic equations in Ito form on manifolds, the Schwartz Morphism \cite{emery} and the Ito Bundle \cite{belopolskaja,gliklikh}, see also \cite{brzezniak}. Focussing specifically on the Schwartz morphism, as studied in Emery (1989) \cite{emery}, Brigo studied its relationship with the 2-jet approach and found them to be very close. The 2-jet interpretation was published in Armstrong and Brigo (2017, 2018) \cite{armstrongjetsgsi,armstrongrspa}, which led next to Armstrong and Brigo investigating how one could project a stochastic differential equation on a sub-manifold in an optimal way. Based on Ito Taylor expansions, two different projections satisfying two different types of optimality were found, the Ito-vector and the Ito-jet projections. The Ito-jet projection is superior in terms of optimality, in that it has a higher order of optimality in a precise sense. These results were presented at ICMS in Edinburgh by Armstrong and Brigo (2016) \cite{armstrongbrigoicms}, at a conference co-organized in 2015 again by Dodson, almost thirty years after the original conference he had organized in Lancaster where the projection filter idea had first been presented. The ICMS conference was co-organized with Frank Critchley and Frank Nielsen. The two projections were studied further and some technical problems concerning tubular neighborhoods were solved with the help of the PhD student Emilio Rossi Ferrucci, leading to the publication Armstrong, Brigo and Rossi Ferrucci (2019) \cite{armstronglms}, see also Armstrong, Brigo and Rossi Ferrucci (2018) \cite{rossiferrucci}, where Rossi Ferrucci helped re-derive the optimal projections through  constrained optimizations as opposed to Ito Taylor expansions, and where ambient coordinates are used. 

In this last paper \cite{armstronglms}, information geometry comes back as an application of the now optimal projections both in Hellinger distance and direct $L^2$ metric, comparing them in a numerical case with the traditional Stratonovich projection of previous works. It turns out that Stratonovich is also optimal for a particular criterion that is, however, not a particularly interesting or natural one, so that the Ito-jet projection filter should be preferred in general.  
 
In this paper we will first present a literature review of projection filtering as done by other authors, following the original papers \cite{brigo1,brigo2}, and then we will detail the above history as much as possible, with derivations, equations and references. We will finally sketch future problems where information geometry might give a contribution. For the reader convenience, we summarize the different projection filtering approaches in Table \ref{tab:pfclass}.  
 
 \begin{table}[h!]
 \begin{tabular}{|c|cc|}
 \hline
Metric $\rightarrow$      & Hellinger & Direct $L^2$     \\
Projection $\downarrow$  & & \\ 
\hline
Stratonovich & Stratonovich ``classic'' PF  & Stratonovich ``classic'' PF \\
projection & exponential families \cite{brigo1,brigo2} & mixture families \cite{armstrongbrigomcss} \\ \hline
Ito-vector &  Optimal Ito-vector PF  & Optimal Ito-vector PF \\
projection &  Gaussian familiy \cite{armstronglms} & Gaussiam family \cite{armstronglms} \\ 
              &  exponential families?  & mixture families? \\ \hline
Ito-jet  & Optimal Ito-jet PF  & Optimal Ito-jet PF  \\
projection &  Gaussian familiy \cite{armstronglms} & Gaussiam family \cite{armstronglms} \\ 
              &  exponential families?  & mixture families? \\ \hline
 \end{tabular}
 \caption{A simplified classification of projection filters (PFs).}\label{tab:pfclass}
 \end{table}

\section{Other works based on the classic projection filters}\label{sec:otherpf} 
 
 Our original work on projection filters was further studied and applied to several fields by subsequent authors. Here we mention only a few examples to illustrate the breadth of the possible use of information geometry and dynamics in applications.  
 
Eagle and Soatto (2011) \cite{soatto} briefly mention the projection filter as one of the possible algorithms for on-line estimation in the context of visual-inertial navigation, mapping and localization.
Lermusiaux (2006) \cite{ocean} mentions the projection filter as a possible tool for estimation of uncertainties for ocean dynamics. 
Kutschireiter, Rast, and Drugowitsch (2022)\cite{circular} apply the projection filter to continuous time circular filtering. 
Projection filters have been applied to quantum systems for example in van Handel and Mabuchi (2005) \cite{quantum1} and in Gao, Zhang and Petersen (2019) \cite{quantum2}. 
Ma, Zhao, Chen and Chang (2015) \cite{hazardposition} apply projection filters to hazard position estimation.
Vellekoop and Clark (2006)\cite{changepoint} extend the projection filter theory to deal with changepoint detection. 
Tronarpand and S\"{a}rkk\"{a}  (2019) \cite{arbitrarylike} present a projection filter for systems with discrete time measurement having arbitrary likelihoods.     
Surace and Pfister (2017) \cite{MLEPartial} apply the Gaussian projection filter to estimate the parameters of a partially observed diffusion. 
Harel, Meir and Opper (2015)\cite{neuralencoding} apply the assumed density filters, equivalent to the projection filters, to the filtering of optimal point processes with applications to neural encoding.  
Azimi-Sadjadi and Krishnaprasad (2005)\cite{navigation} apply projection filter algorithms to navigation. 
Br\"ocker and Parlitz (2000)\cite{chaosfiltering} apply projection filter techniques to address noise reduction in chaotic time series. 
Zhang, Wang, Wu and Xu (2014)\cite{fiber} apply the Gaussian projection filter as part of their estimation technique to deal with measurements of fiber diameters in melt-blown nonwovens.
The projection filter further attracted the attention of the Swedish Defense Research Agency, that summarized and studied it in 2003 in the report \cite{foi}. 

\section{Optimal projection of stochastic differential equations}\label{sec:optprojgen}

To study the optimal approximation of the filtering problem equation, we first introduce three different types of projections on submanifolds for SDEs. For simplicity we illustrate the differences with $n$-dimensional submanifolds $M$ of $\mathbb{R}^r$ but we will then apply this with $L^2$ replacing $\mathbb{R}^r$ and our chosen statistical manifold replacing $M$. 

\begin{figure}[h]
  \begin{minipage}[c]{0.5\textwidth}
    \includegraphics[width=\textwidth]{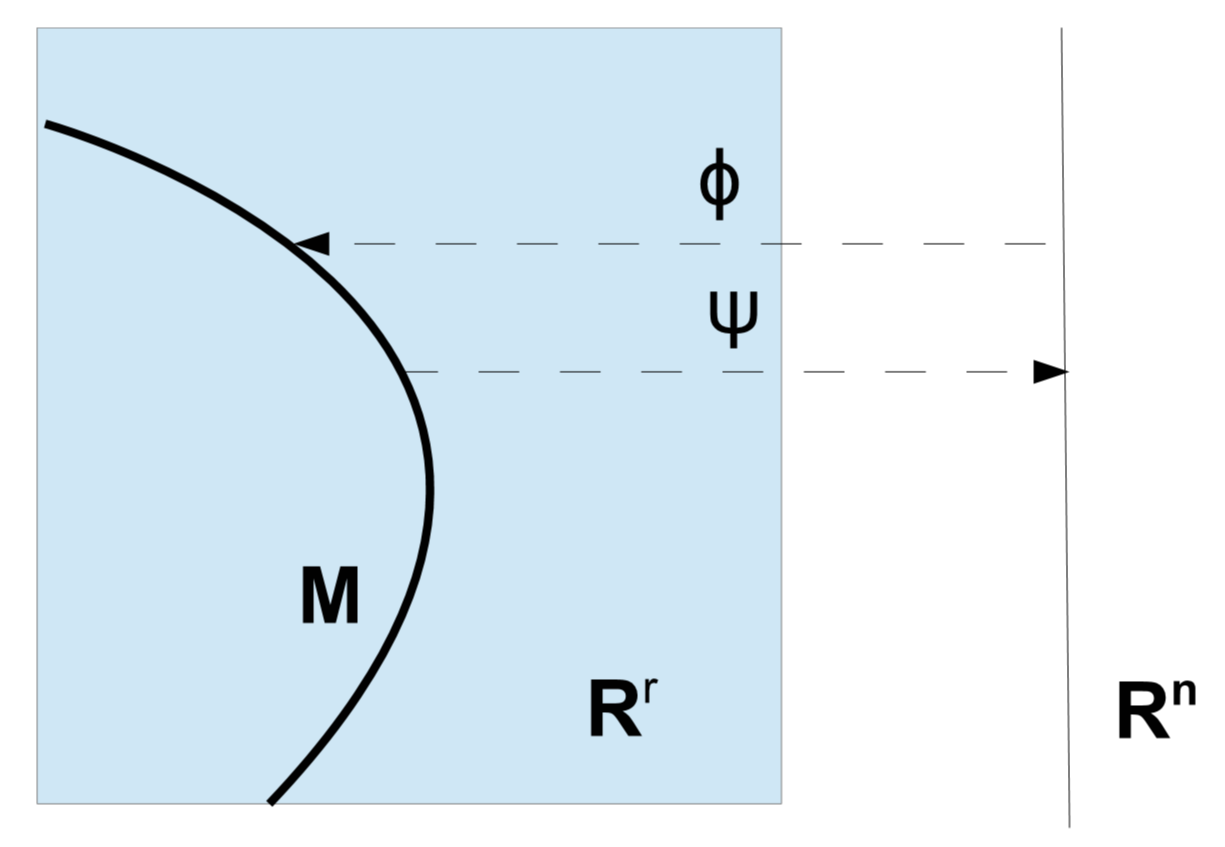}
  \end{minipage}\hfill
  \begin{minipage}[c]{0.4\textwidth}
    \caption{
    Approximate SDE\\ $dX= a(X) dt + b_\alpha(X) dW^\alpha$\\ in $\mathbb{R}^r$
with $\phi(Y) \in M$, where \\ $dY = A(Y) dt + B_\alpha(Y) dW^\alpha$\\ is SDE in $\mathbb{R}^n$, $X_0 = \phi(Y_0) \in M$,  $n$-dimensional manifold of $\mathbb{R}^r$. $\psi$ is the specific chart, $\phi=\psi^{-1}$. 
    } \label{fig:manifold}
  \end{minipage}
\end{figure}

Figure \ref{fig:manifold} illustrates our SDE setting. 

Here we will keep the exposition simple and will not provide all the mathematical details. The full exposition would require tubular neighborhoods, exit times, and constrained optimizations, see \cite{rossiferrucci}. 


Our problem: given a SDE for $X$ on $\R^r$, $dX= a(X) dt + b_\alpha(X) dW^\alpha$ (Einstein convention applies, so $b_\alpha(X) dW^\alpha = \sum_\alpha b_\alpha(X) dW^\alpha$), with $M \subset \R^r$ an $n$-dimensional manifold of $\R^r$, and $X_0 \in M$, we wish to find a SDE $\phi(Y)$ in $M$ starting at $\phi(Y_0)=X_0$,  $dY = A(Y) dt + B_\alpha(Y) dW^\alpha$ being the $\R^n$ coordinates SDE,  whose solution approximates $X$ in some optimal way. Clearly, $r > n$. For brevity, we will often write $a$ for $a(X)$ and $b$ for $b(X)$, and similarly $A$ for $A(Y)$ and $B$ for $B(Y)$. 

In the Stratonovich projection, the above Ito SDE is first transformed in Stratonovich form, obtaining the Stratonovich SDE for $X$ on $\R^r$, 

$dX= \bar{a}(X) dt + b_\alpha(X)\circ dW^\alpha$, where $\bar{a}$ is the transformed drift for the Stratonovich version. 

The full derivation of the optimal approximation and the theory on which this is based is given in full detail in Armstrong, Brigo and Rossi Ferrucci (2019, 2018)\cite{armstronglms, rossiferrucci}, to which we refer for full details. 

\subsection{Stratonovich projection  via tangent space projection}

\begin{figure}[h!]
  \begin{minipage}[c]{0.5\textwidth}
    \includegraphics[scale=0.2]{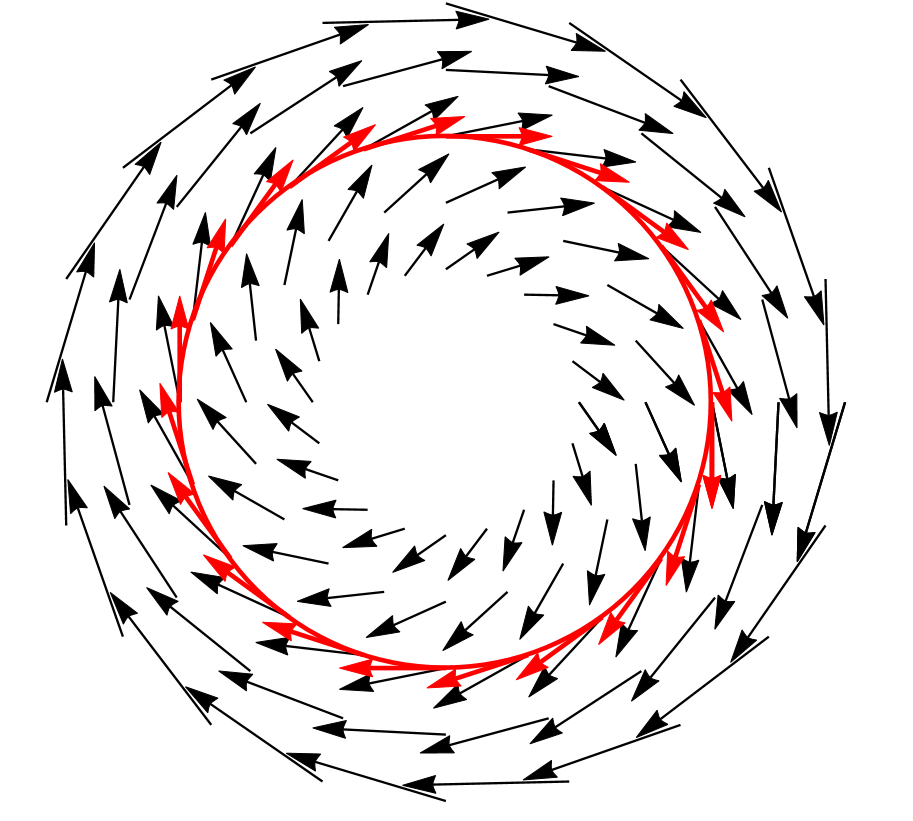}
  \end{minipage}\hfill
  \begin{minipage}[c]{0.5\textwidth}
    \caption{A pictorial representation of the tangent projection $\Pi$ of vector fields on a (red) circle that is used directly in the Stratonovich and Ito-vector projections. 
    } \label{fig:tangentproj}
  \end{minipage}
\end{figure}

\begin{figure}[h!]
  \begin{minipage}[c]{0.5\textwidth}
\includegraphics[scale=0.28]{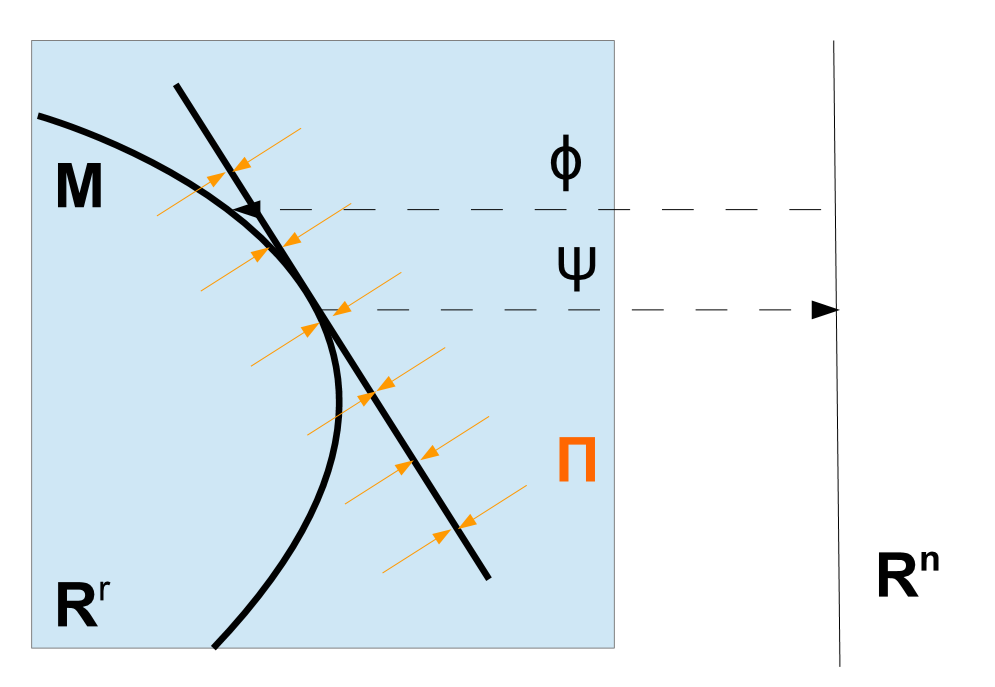}
  \end{minipage}\hfill
  \begin{minipage}[c]{0.5\textwidth}
    \caption{
    {\bf Stratonovich projection}. Write the Ito SDE for $X$ in Stratonovich form \\ 
$dX = \bar{a} \ dt + b_\alpha \circ dW^\alpha \ \ \mbox{in} \ \ \R^r$, \\
 $X_0 \in M$. Apply the tangent space projection $\Pi$ on $M$ to obtain the $M$-SDE\\ for $Z=\phi(Y), Z_0=X_0$,\\
$dZ = \Pi_Z[\bar{a}] \ dt + \Pi_Z[b_\alpha] \circ dW^\alpha$
} \label{fig:stratproj}
  \end{minipage}
\end{figure}

The Stratonovich projection is based on projection on the tangent space of the manifold $M$ of the  vector fields $\bar{a}$ and $b_\alpha$'s. This is illustrated in Figure~\ref{fig:stratproj}. 
For an expression in coordinates we have that the resulting Stratonovich projection SDE is
$dY = \hat{A}(Y) dt + B_\alpha(Y)\circ dW^\alpha$ where
\begin{eqnarray*}
B_\alpha(Y_t,t) &=& (\psi_*)_{\phi(Y_t)} \Pi_{\phi(Y_t)} b_\alpha(\cdot,t),\\
\hat{A}(Y_t,t) &=& (\psi_*)_{\phi(Y_t)} \Pi_{\phi(Y_t)} \bar{a}(\cdot,t)  .
\end{eqnarray*}
As one can see, in the Stratonovich projection we project the $dt$ and $dW_t$ vector fields separately, and get a new Stratonovich SDE that stays on the manifold. This is the projection that was used in the classic projection filter works (1987-2016) we mentioned in Sections~\ref{sec:clastratpr} and \ref{sec:classtratl2}. Also, all the references in Section~\ref{sec:otherpf} consider the Stratonovich projection filters.  The use of the Stratonovich projection as a first method to derive a finite dimensional evolution equation is understandable, given its simplicity. One simply projects the two vector fields separately, and obtains automatically a new SDE on the submanifold. This is the main reason why one uses the Stratonovich projection. It also behaves well under change of coordinates, as coefficients are really vector fields. 

However, given that Stratonovich SDEs do not behave as well as Ito SDEs probabilistically, what does this projection achieve for the SDE solution {\emph{as a whole}}, rather than for the $dt$ and $dW$ vector fields taken separately? How do we justify the Stratonovich projection for $X$ as a whole?

Initially, when struggling with the problem of optimality, we thought the Stratonovich projection had no optimality. A first but somewhat weak  justification for the Stratonovich projection is that for $b=0$ it coincides with the optimal ODE projection, minimizing the leading term of the Taylor expansion for $|\phi(Y)-X|^2$. However, if $b\neq 0$ and we are really dealing with a SDE rather than a ODE, is there any other sense where the Stratonovich projection is optimal? We found out there is a criterion, after all, for which this projection is optimal. 

The Stratonovich projected SDE has a solution minimizing the leading ($t$-term) coefficient of the Taylor expansion of  
\[\frac{1}{2} \left( \E [ | {X_{-t}} -  \phi(Y_{-t})|^2 ]  +
\E [ | {X_t} -  \phi(Y_t)|^2 ]  \right) \]
or of
\[ \frac{1}{2} \left( \E [ |{\pi(X_{-t})}-\phi(Y_{-t})|_r^2 ] +
\E [ |{\pi(X_t)}-\phi(Y_t)|_r^2 ] \right) \]
where $\pi$ is the metric projection from $\R^r$ on $M$, where in general if 
\[ d \xi_t = a(t,\xi_t) dt + b_\alpha(t,\xi_t) \circ dW^\alpha \]
then 
\[ \xi_{-t} := \zeta_t , \ \ d \zeta_t = - a(t,\zeta_t) dt - b_\alpha(t,\zeta_t) \circ d\bar{W}^\alpha\]
where $\bar{W}$ is a new Brownian motion independent of $W$. 
Clearly this is inspired by time symmetry, it is not very natural and is actually quite ad hoc, using time $0$ and the initial condition as an anchor point. This is not very helpful in practical applications. It is a somewhat artificial criterion of optimality that is not very interesting.  

\subsection{Ito-vector projection via tangent space projection}
\begin{figure}[h!]
  \begin{minipage}[c]{0.5\textwidth}
\includegraphics[scale=0.28]{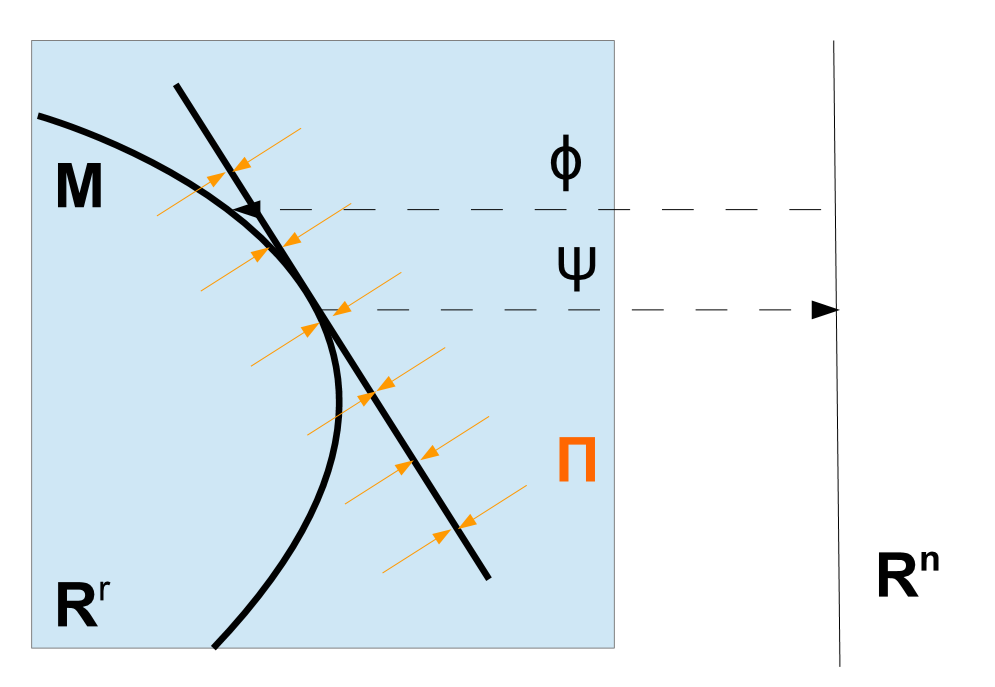}
  \end{minipage}\hfill
  \begin{minipage}[c]{0.5\textwidth}
    \caption{
{\bf It\^o vector projection.} Minimize leading ($t$-term) coeff. of Taylor expansion of  
$\E [ | {X_t} -  \phi(Y_t)|^2 ]$,
to get $B$, {\bf but that coefficient does not vanish}, so the error stays order $t$ rather than $t^2$. Fixing that $B$ in a neighborhood, get $A$ by minimizing the next leading term coeff ($t^2$-term).
This results also in $A$ minimizing (regardless of $B$), up to order $t$,    
the weak error $| \E [ {X_t} -  \phi(Y_t) ]|^2$. 
} \label{fig:itovecproj}
  \end{minipage}
\end{figure}

The idea of the Ito vector projection and the criterion it minimizes are summarized in Figure~\ref{fig:itovecproj}.

The coefficients $A$ and $B$ for the SDE on $\R^n$ such that $\phi(Y)$ is the Ito-vector projection of $X$ on $M$ are 

\begin{equation*}
B_\alpha(Y_t,t) = (\psi_*)_{\phi(Y_t)} \Pi_{\phi(Y_t)} b_\alpha(\cdot,t)\ \  \mbox{(same as Stratonovich proj.)}
\end{equation*}
\[
A(Y_t,t) = (\psi_*)_{\phi(Y_t)} \Pi_{\phi(Y_t)} \left(a(\cdot,t) - \frac{1}{2}(\D_{B_\alpha(Y_t,t)} \phi_*) B_\beta(Y_t,t) g^{\alpha\beta}_E \right).
\]
The term $B$ is essentially the same as in the Stratonovich projection, consisting of the projection of the vector field $b$ on the tangent space of the manifold $M$. 

The matrix $g^{\alpha,\beta}_E = 1_{\{\alpha = \beta\}}$, or more generally it is the symmetric $2$-form defining the Euclidean metric on $\R^m$ in the non-orthonormal case. It may be helpful to recall that in Euclidean space notation (with orthonormal coordinates) we have 
\[ (\nabla_{B_\alpha} f_{i,\ast} ) B_\beta g^{\alpha,\beta}_E = \mbox{Trace}[ B^T (H f_i) B] \]
with $H$ the Hessian operator. Given the drawback 	pointed out in Figure~\ref{fig:itovecproj} of this projection being optimal only up to order $t$ rather than $t^2$, can we find another projection that, differently from the It\^o vector projection, is {\emph{consistently optimal}} up to order $t^2$?

\begin{figure}[h!]
  \begin{minipage}[c]{0.5\textwidth}
    \includegraphics[scale=0.3]{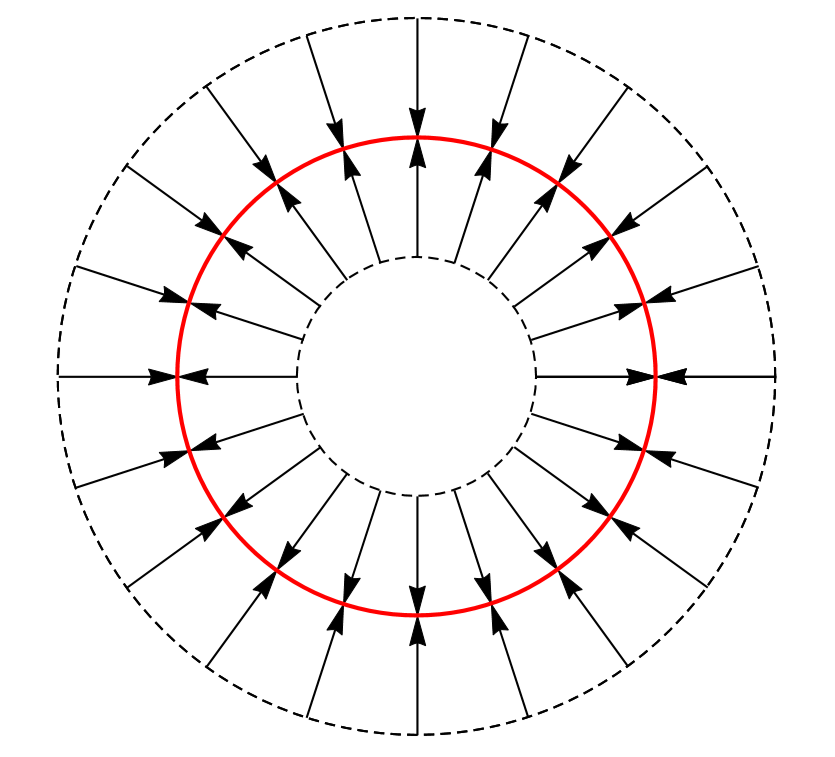}
        \includegraphics[scale=0.3]{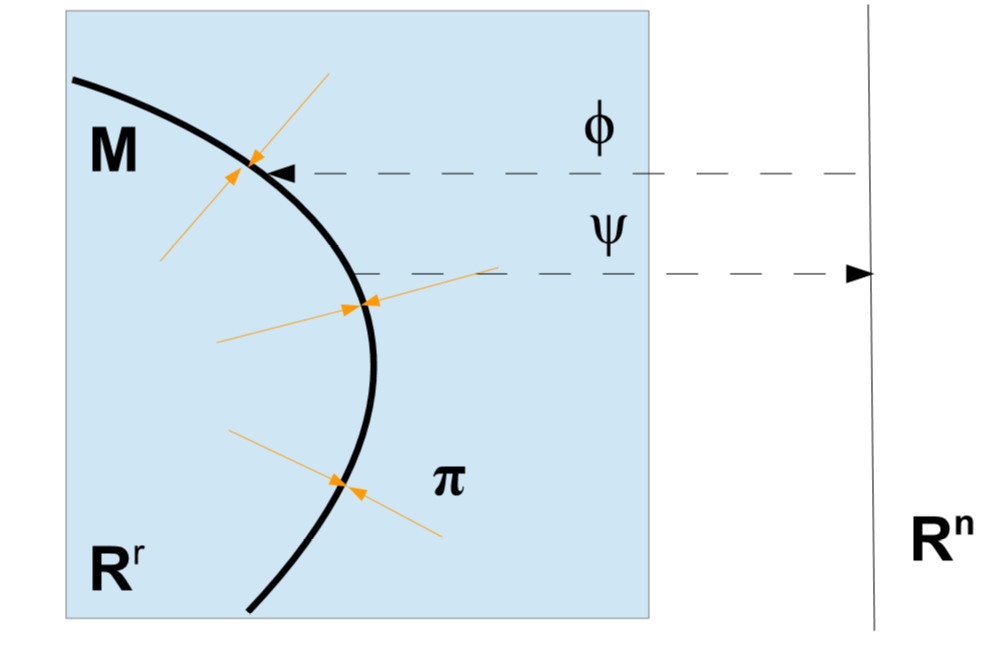}
  \end{minipage}\hfill
  \begin{minipage}[c]{0.45\textwidth}
    \caption{Above: A pictorial representation of the metric projection $\pi$ on a (red) circle $M$ that is used in the Ito-jet projection. \\ Below: $\pi$ is the metric projection\\ $\pi: \R^r \rightarrow M$,\\  defined on a tubular neighborhood of $M$, of which the  earlier linear projection  $\Pi$ is the first order  component.\\ Using the chart $\psi$, set $\tilde{\pi} = \psi\circ \pi$.
In the {\bf It\^o jet projection} we make the $t$ coefficient of the expansion {\emph{vanish}} and minimize the leading $t^2$ coefficient of the Taylor expansion for the error
$\E [ d_M({\pi(X_t)},\phi(Y_t))^2 ]$ or \\ 
$E [ |{\pi(X_t)}-\phi(Y_t)|_r^2 ]$  for small $t$, 
so that we attain optimality up to order $2$. $d_M$ is the distance on $M$.} \label{fig:metricproj}
  \end{minipage}
\end{figure}

The third and most optimal of the three projections will achieve this and will be based on the metric projection, rather than the tangent projection. A pictorial representation of the metric projection is presented in Figure~\ref{fig:metricproj}.

\subsection{Ito jet projection via metric projection}\label{sec:itojetp}

With this criterion, $B$ is as before, as in the Stratonovich and Ito vector projections, so we do not repeat it, while $A$ is given by  
\[
A(Y_t,t)= \tilde{\pi}_*(a(\phi(Y_t),t))
+ \frac{1}{2} (\D_{b_{\alpha(\phi(Y_t),t)}}\tilde{\pi}_*)b_\beta(\phi(Y_t),t) g^{\alpha\beta}_E.
\]
In Armstrong, Brigo and Ferrucci-Rossi (2019) \cite{armstronglms} we show the detailed calculations for $A$, which are very long and laborious. They involve the expansion of the metric projection. In short, define the metric tensor in a tubular neighborhood of $M$, 
\begin{equation}
h_{a b} = \frac{\partial \phi^\alpha}{\partial x^a}
    \frac{\partial \phi^\alpha}{\partial x^b}
    \label{metricTensor}
\end{equation}    
and call $h^{a b}$ the inverse of $h$. 
The differential ${\tilde{\pi}}_*$ of ${\tilde{\pi}}$ is well known to be given
by the linear projection onto $\Image \phi_*$ composed with the map $\phi_*^{-1}$. Hence
${\tilde \pi}_*$ is the unique linear map with ${\tilde \pi}_* \circ \phi_*$ equal to the identity and with
kernel equal to the orthogonal complement of $\Image \phi_*$. We deduce that ${\tilde \pi}_*$ has
the following components:

\begin{equation}
 \Pi^a_b:= ({\tilde \pi}_*)^a_b =  \frac{\partial \phi^b}{\partial x^\alpha} h^{a \alpha}, \qquad a \leq n, \alpha \leq n, b \leq r.
\label{piFirstOrder}
\end{equation}

We note that the differential or tangent map ${\tilde \pi}_*$ is the best linear approximation of the metric projection ${\tilde \pi}$ around the relevant point $x = \phi(y) \in M$, and it coincides with the classic linear projection $\Pi_{\phi(y)}$ on the tangent space of $M$.
%

\[
A^i = \Pi^i_\alpha a^\alpha + 
\left( -\frac{1}{2} \frac{ \partial^2 \phi^\gamma}{\partial x^\alpha \partial  x^\beta}
\Pi^i_\gamma
\Pi^\alpha_\delta  
\Pi^\beta_\epsilon
+ \frac{ \partial^2 \phi^\epsilon}{ \partial x^\alpha \partial x^\beta}
\Pi^\beta_\delta
h^{i \alpha} 
-  \frac{ \partial^2 \phi^\gamma}{ \partial x^\alpha \partial x^\beta}
\Pi^\beta_\epsilon
\Pi^\eta_\gamma
\Pi^\zeta_\delta
h_{\eta \zeta} h^{i \alpha}
\right)  \times b^\delta_\kappa b^\epsilon_\iota [W^\kappa, W^\iota]_t.
\]

We summarize the different projections and the optimality criteria used to determine their drifts in \Cref{tab:projtypes} in the conclusions. The diffusion coefficient is identical for all three projections.

\section{Application to Non-linear Filtering via Information Geometry}
\label{applicationSection}
We studied the application of the new projections to nonlinear filtering via information geometry in Armstrong, Brigo and Rossi Ferrucci (2019) \cite{armstronglms}. Here we summarize the results of that paper, showing how our new projection methods work for stochastic filtering. As explained in the introduction, this enhances optimality of the approximations compared to our previous works in \cite{brigo1}, \cite{brigo2} and \cite{armstrongbrigomcss}.

Let us first summarize the filtering problem for diffusions. One has a signal $X$ that evolves according to a SDE, and observes a process $Y$ which is a function of this signal plus noise. This is standard notation, but these $X$ and $Y$ are not to be confused with the processes we used earlier in the paper, in that they are not the $\R^r$ process to be approximated and its $\R^n$ approximation.

The filtering problem consists in estimating the signal $X$ given the present and past observations $Y$. If $t$ is the current time, the solution of the filtering problem is the probability density of the state $X_t$ conditional on the observations from time 0 to time $t$, call this density $p_t$. The density $p_t$ follows the Kushner-Stratonovich (or alternatively the Zakai) stochastic partial differential equation (SPDE) that, under some technical assumptions, can be seen as a stochastic differential equation in the infinite dimensional $L^2$ space of square roots of densities (Hellinger metric) or of densities themselves (direct $L^2$ metric).

The process we wish to approximate on a low dimensional manifold is $p_t$, which represents the $X_t$ of our earlier sections. The $\R^r$ space of our earlier sections is the $L^2$ infinite dimensional space, while the submanifold $M$ is a finite dimensional family of probability densities parametrized by $\theta$, acting as coordinates: $\{p(\cdot,\theta), \ \theta \in \Theta \subset \R^n\}$. $\theta_t$ plays the role of what we were calling $Y_t$ earlier in the paper. For brevity, we will often omit the argument and write $p(\cdot,\theta)=p_\theta$.  
We aim at finding a SDE for $\theta$ such that $p_{\theta_t}$ approximates the optimal filter $p_t(\cdot)$ in an optimal way. Note that in the previous part of the paper we had a dimensionality reduction from $r$ to $n$, whereas now we go from infinite dimensional $p_t$ to $n$-dimensional $\theta_t$.

We discuss in \cite{armstronglms} how the setting can be generalized from $\R^r$ to $L^2$.

\subsection{The Kushner Stratonovich equation}
\label{filteringSubsection}

We suppose that the state $X_t \in \R^m$ of a system
evolves according to the equation:
\[ \ed X_t = f(X_t,t) \, \ed t + \sigma(X_t,t) \, \ed W_t \]
where $f$ and $\sigma$ are smooth $\R^m$ valued functions
and $W_t$ is a Brownian motion. One typically adds growth conditions to ensure a global existence and uniqueness result for the signal equation, see for example \cite{armstrongbrigomcss} and references therein for the details.

We suppose that an associated process, the observation process, $Y_t \in \R^d$
evolves according to the equation:
\[ \ed Y_t = b(X_t,t) \, \ed t + \ed V_t \]
where $b$ is a smooth $\R^d$ valued function and $V_t$ is a Brownian motion independent of $W_t$.
Note that the filtering problem is often formulated with an additional constant in terms
of the observation noise. For simplicity we have assumed that the system is scaled so that this can be omitted.

The filtering problem is to compute the conditional
distribution of $X_t$ given a prior distribution for $X_0$
and the values of $Y$ for all times up to and including~$t$.

Subject to various bounds on the growth of the
coefficients of this equation, the assumption that
the distribution has a density $p_t$ and suitable
bounds on the growth of $p_t$ one can show that
$p_t$ satisfies the Kushner--Stratonovich SPDE:
\begin{equation}
\ed p = {\cal L^*}  p \  \ed t 
+ p[b - E_p(b)]^T [ \ed Y - E_p(b) \ed t]
\label{KSIto}
\end{equation}
where $E_p$ denotes the expectation with respect to
the density $p$, 

$E_p[f] = \int f(x) p(x) dx$, and the forward diffusion operator ${\cal L}^*_t$ is defined by:  
\begin{equation}
{\cal L}_t^* \phi = - \frac{\partial}{\partial x^i} [ f_i(x,t) \phi ] + \frac{1}{2} \frac{\partial^2}{\partial x^i \partial x^j} [a_{ij}(x,t)\phi]
\end{equation}
where $a=\sigma \sigma^T$.
Note that we are using the Einstein summation convention in this expression. 

In the event that the coefficient functions $f$
and $b$ are all linear and $\sigma$ is a deterministic function of time 
one can show that
so long as the prior distribution for $X$ is Gaussian, or deterministic, 
the density $p$ will be Gaussian at all subsequent times.
This allows one to reduce the infinite dimensional 
equation \eqref{KSIto} to a finite dimensional stochastic differential equation for the mean and covariance matrix of this normal distribution. This finite dimensional
problem solution is known as the Kalman filter.

For more general coefficient functions, however, equation
\eqref{KSIto} cannot be reduced to a finite dimensional
problem \cite{hazewinkel}. Instead one might seek approximate solutions of \eqref{KSIto} that belong to some given statistical family of densities. This is
a very general setup and includes, for example,
approximating the density using piecewise linear
functions to derive a finite difference approximation
or approximating the density with Hermite polynomials
to derive a spectral method. Other examples include exponential families (considered in \cite{brigo2,brigo1}) and mixture families (considered in \cite{armstrongBrigo,armstrongbrigomcss}).

Our projection theory tells us how one can find good approximations 
on a given statistical family with respect to a given
metric on the space of distributions. We illustrate
this by writing down the It\^o-vector and It\^o-jet projection of
\eqref{KSIto} for the $L^2$ and Hellinger metrics
onto a general manifold\footnote{Note that it is also possible to consider projecting the Zakai equation. However, as explained in 
\cite{armstrongbrigomcss}, one expects that projecting the Kushner--Stratonovich equation will lead to smaller error terms in direct metric, whereas the projected equations are the same in Hellinger metric. See \cite{armstrongbrigomcss} for a discussion}. 

We will then examine some numerical results regarding 
the very specific case
of seeking approximate solutions using Gaussian distributions. The idea of approximating the solution
to the filtering problem using a Gaussian
distribution has been considered by numerous authors
who have derived variously, the extended Kalman filter
\cite{pardoux}, assumed density filters \cite{kushner} and Stratonovich projection filters \cite{brigo1}. Some of these are related, for example the assumed density filters and Stratonovich projection filters in Hellinger metrics for Gaussian (and more generally exponential) families coincide \cite{brigo2}.  Using our new projection methods,  we will be able to derive projection filters which
 outperform all these other filters
(assuming performance is measured over small time intervals using the appropriate Hilbert space metric).

We note that \eqref{KSIto} is an infinite dimensional SDE
driven by a continuous semi-martingale. The definitions and results
given in Section \ref{sec:optprojgen} were only stated in the
finite dimensional case for SDEs driven by Brownian motion.
The definition of It\^o--Taylor series can be generalized straightforwardly
to this situation and hence the definition of the It\^o projections can
be applied in this context also.

More generally, for the the geometry of approximations to the infinite dimensional filtering problems based on $L^2$ or Orlicz charts we refer for example to \cite{brigo1,brigo2,armstrongbrigomcss,armstronglms,brigopistone,brigopistone2,newton1,newton2}.

\subsection{Stratonovich projections}
The Stratonovich projection filters have been abundantly studied in \cite{brigo1,brigo2} in Hellinger metric, and in  \cite{armstrongbrigomcss} in direct metric, see also references in Section~\ref{sec:otherpf} for the Hellinger case. Here we briefly summarize them. 
For this method, the optimal filter SPDE is given by putting the optimal filter equation \eqref{KSIto} in Stratonovich form, obtaining 
\begin{equation}\label{KSstrat}
   dp = {\cal L}^\ast\, p\,dt
   - \frac{1}{2}\, p\, [\vert b \vert^2 - E_{p}\{\vert b \vert^2\}] \,dt
   +  p\, [b-E_{p}\{b\}]^T \circ dY\ .
\end{equation}

\subsubsection{Stratonovich projection filter in $L^2$ direct metric}

We use first the $L^2$ direct metric, by assuming densities are squared integrable, and view densities are elements of $L^2$ with the related inner product.
For a discussion on conditions under which a unnormalized version of the density $p_t$ of the optimal filter \eqref{KSIto} is in $L^2$ (Zakai Equation) see for example \cite{ahmedbook}. 

In the case of two probability
density functions $p(x)$ and $q(x)$ on $\R^m$, we assume them to be in $L^2$, and the $L^2$ squared direct distance is
\[ 
 \int ({p(x)} - {q(x)})^2 \, \ed x .
\]
We wish to consider an $n$-dimensional family
of distributions $p_\theta$ parameterized by $n$ real valued parameters $\theta^1$, $\theta^2$, $\ldots$, $\theta^n$. For example, for a scalar state space $x$, one may consider the 2 dimensional Gaussian family:
\begin{equation}
 p_\theta(x) = \frac{1}{(\theta^2) \sqrt{2 \pi}}\exp\left(-\frac{(x-(\theta^1))^2}{2 (\theta^2)^2}\right). 
 \label{gaussianFamily}
\end{equation}
For more general exponential or mixture families see \cite{brigo1, brigo2, armstrongbrigomcss}. 
Note that we have chosen to follow the differential
geometry convention and use upper indices for
the coordinate functions $\theta^i$ so we have 
been careful to distinguish powers from indices using
brackets.

Once the family is chosen, in the Stratonovich projection we can project directly the $dt$ and the $dY$ terms of the optimal filter equation \eqref{KSstrat} on the tangent space of the chosen manifold $p_\theta$ in direct metric, as both terms behave as vector fields. 

More formally, an $n$-dimensional family is given
by a smooth embedding $\phi:\R^n \rightarrow L^2(\R^m)$.
The tangent vectors $\phi_*\frac{\partial}{\partial \theta^i} \in L^2(\R^m)$ are simply the partial derivatives
\[
\frac{\partial p_\theta}{\partial \theta^i}.
\]

Let us write:
\[ g_{ij}(\theta)=\int_{\R^m} \frac{\partial p_\theta(x)}{\partial \theta^i}
\frac{\partial p_\theta(x)}{\partial \theta^j} \, \ed x.\]
This defines the induced metric tensor on the manifold
$\phi(\R^n)$. We will write $g^{ij}$ for the inverse of the matrix $g_{ij}$. The projection operator $\Pi_{\phi(\theta)}$
is then given by
\begin{equation*}
\begin{split}
\Pi_{\phi(\theta)} (v) &= 
\sum_{i,j=1}^n g^{ij}(\theta) \left\langle v, \phi_* \frac{\partial}{\partial \theta^i} \right\rangle_{L^2}
\phi_*\frac{\partial}{\partial \theta^j} \\
&= \sum_{i,j=1}^n g^{ij}(\theta) \left( \int_{\R^m} v(x) \frac{\partial p_\theta(x)}{\partial \theta^i} \, \ed x \right)
\phi_*\frac{\partial}{\partial \theta^j} \, .
\end{split}
\end{equation*}
Thus
\[
\phi_*^{-1} \Pi_{\phi(\theta)} (v)
= \sum_{i,j=1}^n g^{ij}(\theta) \left( \int_{\R^m} v(x) \frac{\partial p_\theta}{\partial \theta^i} \, \ed x  \right)
\frac{\partial}{\partial \theta^j} \, .
\]

We can now write down the Stratonovich projection of \eqref{KSstrat} with respect to the $L^2$ metric. It is:
\[ \ed \theta^i = \hat{A}^i(\theta) \, \ed t + B^i(\theta) \circ \, \ed Y_t \]
where:
\[ B^i(\theta) =
\sum_{j=1}^n g^{ij}(\theta) \left( \int_{\R^m}
(p_\theta(x)(b(x)-E_{p_{\theta}}(b)))^T
\frac{\partial p_\theta(x)}{\partial \theta^j} \, \ed x \right)
\]
and
\begin{small}
\[ \hat{A}^i(\theta) =
\sum_{j=1}^n g^{ij}(\theta) \left( \int_{\R^m}
\left(
{\cal L}^*p_\theta - \frac{1}{2} p_\theta (|b|^2-E_{p_\theta}[|b|^2])
\right)
\frac{\partial p_\theta}{\partial \theta^j} \, \ed x \right).
\]
\end{small}

\subsubsection{Stratonovich projection filter in Hellinger metric}

If we use instead projection under the Hellinger distance, we need to work with square roots of densities. 
Indeed, the Hellinger metric is a metric on probability
measures. In the case of two probability
density functions $p(x)$ and $q(x)$ on $\R^m$, that now need only be in $L^1$, the Hellinger distance is
given by the square root of:
\[ 
\frac{1}{2} \int (\sqrt{p(x)} - \sqrt{q(x)})^2 \, \ed x .
\]
In other words, up to the constant factor of $\frac{1}{2}$
the Hellinger metric corresponds to the $L^2$ norm on
the square root of the density function rather than on the density itself (as in the previous case of the direct metric). The Hellinger metric has the important advantage of making the metric independent of the particular background density that is used to express measures as densities. The $L^2$ direct distance introduced earlier does not satisfy this background independence.

To get the equation for the square root of the optimal filter density $p$, we can now use the chain rule, as Stratonovich calculus satisfies it. We get
\[ d \sqrt{p} = \frac{1}{2 \sqrt{p}} \circ \ dp \] 
or
\begin{equation}\label{eq:StratKSE} d \sqrt{p} = \frac{1}{2\sqrt{p}}{\cal L}^\ast\, p\,dt
   - \frac{1}{4 }\, \sqrt{p}\, [\vert b \vert^2 - E_{p}\{\vert b \vert^2\}] \,dt
   + \frac{1}{2}\sqrt{ p}\, [b-E_{p}\{b\}]^T \circ dY\ . 
\end{equation}
A family of distributions now corresponds to an embedding
$\phi$ from $\R^n$ to $L^2(\R^m)$ but now $p_\theta=\phi(\theta)^2$. The tangent space is spanned by
the vectors:
\[ \phi_* \frac{\partial}{\partial \theta^i}
 = \frac{ \partial \sqrt{p_\theta} }{\partial \theta^i}. \]
We define a metric on the tangent space by:
\[
h_{ij}(\theta) = \int_{\R^m} \frac{ \partial \sqrt{p_\theta(x)}}{\partial \theta^i} \frac{ \partial \sqrt{p_\theta(x)} }{\partial \theta^j} \, \ed x.
\]
We write $h^{ij}$ for the inverse matrix of $h_{ij}$. The projection operator with respect to the Hellinger metric is:
\begin{equation*}
\begin{split}
\Pi_{\phi(\theta)} (v) 
&= \sum_{i,j=1}^n h^{ij}(\theta) \left( \int_{\R^m} v(x) \frac{\partial \sqrt{p_\theta(x)}}{\partial \theta^i} \, \ed x \right)
\phi_*\frac{\partial}{\partial \theta^j} .
\end{split}
\end{equation*}
We can now write down the Stratonovich projection of \eqref{eq:StratKSE} with respect to the Hellinger metric. It is: 
\[ \ed \theta^i = \hat{A}^i(\theta) \, \ed t + B^i(\theta)\, \circ \ed Y_t \]
where:
\[ B^i(\theta) =
\sum_{j=1}^n h^{ij}(\theta) \left( \int_{\R^m}
\frac{1}{2} \sqrt{p_\theta(x)}(b(x) - E_{p_\theta}(b))^T
\frac{\partial \sqrt{p_\theta(x)}}{\partial \theta^j} \, \ed x. \right)
\]
and
\[
\hat{A}^i(\theta) =
\sum_{j=1}^n h^{ij}(\theta) \left( \int_{\R^m}
\left(
\frac{{\cal L}^* p_\theta(x)}{ 2 \sqrt{p_\theta(x)}}- \frac{1}{4 }\, \sqrt{p_\theta(x)}\, [\vert b(x) \vert^2 - E_{p_\theta}\{\vert b \vert^2\}]\right)
\frac{\partial \sqrt{p_\theta(x)}}{\partial \theta^j} \, \ed x. \right).
\]
This is the Stratonovich SDE for $\theta_t$, the Stratonovich projection of the optimal filter on our family $p_\theta$ in Hellinger metric.

In the next sections we will derive the optimal Ito projections. These will be somehow more complicated and we will thus omit functions arguments such as $\theta$ and $x$ where obvious.

\subsection{It\^o-vector projections}
\label{vectorCalculations}
\subsubsection{The It\^o-vector projection filter in the $L^2$ direct metric}

Let us suppose again that the optimal filter density $p$ lies in $L^2$ and so
we can use the $L^2$ norm to measure the accuracy of
an approximate solution to equation \eqref{KSIto}.

We can now write down the It\^o-vector projection of \eqref{KSIto} with respect to the $L^2$ metric. It is:
\[ \ed \theta^i = A^i \, \ed t + B^i \, \ed Y_t \]
where:
\[ B^i =
\sum_{j=1}^n g^{ij} \left( \int_{\R^m}
(p(b-E_{p(\theta)}(b)))^T
\frac{\partial p}{\partial \theta^j} \, \ed x \right)
\]
and
\begin{small}
\[ A^i =
\sum_{j=1}^n g^{ij} \left( \int_{\R^m}
\left(
{\cal L}^*p - p(b-E_{p(\theta)}(b))^T E_{p(\theta)}(b)
- \frac{1}{2} \sum_{k=1}^n
\frac{\partial^2 p}{\partial \theta^j \partial \theta^k} B^k
\right)
\frac{\partial p}{\partial \theta^j} \, \ed x \right).
\]
\end{small}

\subsubsection{The It\^o-vector projection filter in the Hellinger metric}

Now, to compute the It\^o-vector projection with respect to
the Hellinger metric we first want to write down
an It\^o equation for the evolution on $\sqrt{p}$.

Applying It\^o's lemma to equation \eqref{KSIto} we
formally obtain:
\begin{equation*}
\begin{split}
\ed \sqrt{p}
&= \left( 
\frac{ {\cal L}^* p - p(b - E_p(b))^T E_p(b)}{ 2 \sqrt{p}}
- \frac{p^2 (b - E_p(b))^T(b - E_p(b))}{ 8 p \sqrt{p}}
\right) \, \ed t \\
&{} + \left(
\frac{p(b - E_p(b))^T }{ 2 \sqrt{p}}
\right) \ed Y_t. \\
&= \left( 
\frac{ {\cal L}^* p}{ 2 \sqrt{p}}
- \frac{1}{8}\sqrt{p} (b - E_p(b))^T( b + 3 E_p(b))
\right) \, \ed t \\
&{} + \left(
\frac{1}{2} \sqrt{p}(b - E_p(b))^T
\right) \ed Y_t.
\end{split}
\end{equation*}

A family of distributions now corresponds to an embedding
$\phi$ from $\R^n$ to $L^2(\R^m)$ but now $p=\phi(\theta)^2$. The tangent space is spanned by
the vectors:
\[ \phi_* \frac{\partial}{\partial \theta^i}
 = \frac{ \partial \sqrt{p} }{\partial \theta^i}. \]
We define a metric on the tangent space by:

\[
h_{ij} = \int_{\R^m} \frac{ \partial \sqrt{p} }{\partial \theta^i} \frac{ \partial \sqrt{p} }{\partial \theta^j} \, \ed x.
\]
We write $h^{ij}$ for the inverse matrix of $h_{ij}$. The projection operator with respect to the Hellinger metric is:
\begin{equation*}
\begin{split}
\Pi_{\phi(\theta)} (v) 
&= \sum_{i,j=1}^n h^{ij} \left( \int_{\R^m} v(x) \frac{\partial \sqrt{p}}{\partial \theta^i} \, \ed x \right)
\phi_*\frac{\partial}{\partial \theta^j} .
\end{split}
\end{equation*}
We can now write down the It\^o-vector projection of \eqref{KSIto} with respect to the Hellinger metric. It is: 
\[ \ed \theta^i = A^i \, \ed t + B^i \, \ed Y_t \]
where:
\[ B^i =
\sum_{j=1}^n h^{ij} \left( \int_{\R^m}
\frac{1}{2} \sqrt{p}(b - E_{p(\theta)}(b))^T
\frac{\partial \sqrt{p}}{\partial \theta^j} \, \ed x. \right)
\]
and
\[
\begin{split}
A^i &=
\sum_{j=1}^n h^{ij} \left( \int_{\R^m}
\left(
\frac{{\cal L}^* p}{ 2 \sqrt{p}}
- \frac{1}{8}\sqrt{p} (b - E_{p(\theta)}(b))^T( b + 3 E_{p(\theta)}(b)) \right. \right. \\
&\qquad \left. \left.
- \frac{1}{2} \sum_{k=1}^n
\frac{\partial^2 \sqrt{p}}{\partial \theta^j \partial \theta^k} B^k
\right)
\frac{\partial \sqrt{p}}{\partial \theta^j} \, \ed x. \right).
\end{split}
\]

\subsection{It\^o-jet projections}

Using the formulae from \Cref{sec:itojetp} together with
formulae and techniques analogous to those of \Cref{vectorCalculations} we can
explicitly calculate the It\^o-jet projections of the filtering
equation in both the $L^2$ and Hellinger metrics.

To minimize notation, let us concentrate on the $1$-dimensional state space
filtering problem and project using the $L^2$ metric.

We can formally write the filtering equation in the form:
\begin{equation}
\ed p_t = \mu(p_t) \ed t + \Sigma(p_t) \ed Y_t
\label{infDimEquation}
\end{equation}
where $p_t$ is an $L^2$ function and
\begin{equation}
\label{filteringCoefficientsInfDim}
\begin{split}
\mu(p)(x)&:= \frac{1}{2} \frac{ \ed^2 (\sigma(x)^2 p(x))}{ \ed x^2}
		   - \frac{\ed (f(x) p(x))}{\ed  x} \\
&{}\quad   -  p(x) \left( b(x) - \int_\R p(t) b(t) \ed t \right) \int_{\R} p(t)b(t) \ed t, \\
\Sigma(p)(x)&:= p(x) \left( b(x) - \int_{\R} p(t) b(t) \ed t \right).
\end{split}
\end{equation}
We now suppose that $p_t$ is parameterized as $p_t(x) = \phi(\theta)(x)$
as in \Cref{vectorCalculations}. Using results in  \Cref{sec:itojetp}
we can write down the It\^o-jet projection which is an SDE for the components of $\theta$.

To write down the result it will be useful to define functions $\pi^i(\theta)$ by:
\[
\pi^i(\theta) = h^{ij} \frac{\partial \phi}{\partial \theta^j} (\theta).
\]
We will also use 
angle brackets to denote the $L^2$ inner product. With this
understood, the It\^o-jet projection of the filtering
equations in the $L^2$ metric is given by:
\[ 
\ed \theta^i_t = A^i(\theta) \ed t + B^i(\theta) \ed Y_t
\]
where we have in turn
\[
B^i(\theta) = \langle \pi^i(\theta), \Sigma(\phi(\theta)) \rangle
\]
and
\begin{equation*}
\begin{split}
A^i(\theta) &= \langle \pi^i(\theta), \mu(\phi(\theta)) \rangle \\
& 
\quad - \frac{1}{2}
\left\langle \frac{\partial^2 \phi}{\partial \theta^\alpha \partial \theta^\beta}(\theta),
	      \pi^{i}(\theta)
\right\rangle
\langle
\Sigma(\theta), \pi^\alpha(\theta)
\rangle
\langle
\Sigma(\theta), \pi^\beta(\theta)
\rangle \\
& 
\quad + 
\left\langle \frac{\partial^2 \phi}{\partial \theta^\alpha \partial \theta^\beta}(\theta), \Sigma(\phi(\theta)) \right\rangle
\langle \pi^\beta(\theta), \Sigma(\phi(\theta)) \rangle
h^{i \alpha}(\theta) \\
& 
\quad -
\left\langle \frac{\partial^2 \phi}{\partial \theta^\alpha \partial \theta^\beta}, \pi^\eta(\theta) \right\rangle
\langle \pi^\beta(\theta), \Sigma(\phi(\theta)) \rangle
\langle \pi^\xi(\theta), \Sigma(\phi(\theta)) \rangle
h_{\eta \xi}(\theta) h^{i \alpha}(\theta)
\end{split}
\end{equation*}

The It\^o-jet projection of the filtering equation in the Hellinger metric can be computed in the same way. Indeed we can formally write the filtering equation in the form:
\begin{equation}
\ed q_t = \mu(q_t) \ed t + \Sigma(q_t) \ed Y_t
\end{equation}
where $q_t$ is the square root of the density and the coefficients now
satisfy
\begin{equation}
\label{filteringCoefficientsInfDimHellinger}
\begin{split}
\mu(q)(x)&:= \frac{1}{2 q(x)}
\left(\frac{1}{2} \frac{d^2 (\sigma(x)^2 q(x)^2)}{\ed x^2} - \frac{\ed (f(x)q(x)^2)}{\ed x} \right) \\
&\quad
  - 
 \frac{1}{8} q(x) \left( b(x) - \int_\R q(t)^2 b(t) \ed t \right)
\left( b(x) + 3 \int_{\R} q(t)^2 b(t) \ed t \right), \\
\Sigma(q)(x)&:= \frac{1}{2} q(x) \left( b(x)- \int_{\R} q(t)^2 b(t) \ed t \right).
\end{split}
\end{equation}
Thus we can use the same formulae as above to compute the Hellinger projection
except we must use the coefficients from \eqref{filteringCoefficientsInfDimHellinger} rather than those
from \eqref{filteringCoefficientsInfDim}.

\subsection{Comparison of filters}

In \cite{armstronglms} we compare the different projection filters with each other in a case of cubic sensor perturbing a linear system (where, without perturbation, the Kalman filter would work well). In other words, the state equation is trivial, $dX=dW$, while the observation function is $b(x) = x + \varepsilon x^3$. For small $\varepsilon$, this will be close to a linear system and the extended Kalman filter and other Gaussian filters are supposed to perform well. We make the comparison in \cite{armstronglms}, comparing the different projection filters with the extended Kalman filter and with the Ito assumed density filter (ADF) with assumed Gaussian density.   

First of all, our explicit calculations for the cubic sensor show that the two It\^o projections give rise to new, distinct, Gaussian approximations which are different from all previous filters, so we have indeed introduced new approximate filters. 

The explicit calculations in \cite{armstronglms} show that all the resulting filters for the cubic sensor $b(x,t)=x + \varepsilon x^3$ are equal when $\varepsilon=0$. This provides a basic
sanity check that our formulae correspond to the Kalman filter in
the case of a linear sensor. In general, if we know that the solution
lies in a particular manifold and we project onto that manifold, the
three projection methods will all be exact.

We simulated the example problem $b(x)=x + \varepsilon x^3$ for all of the above
approximate filters with $\epsilon=0.05$. We also computed an ``exact'' solution
using a finite difference method. We define the $L^2$ residual
to be the $L^2$ direct distance between the approximate solution and the ``exact'' solution.
We define the Hellinger residual similarly, as the $L^2$ distance between the square roots of the solution densities.

In \cite{armstronglms} we compare first the $L^2$ residuals for the various methods.
All the projection methods shown are taken using the $L^2$ metric in this case.
The It\^o-vector projection in the $L^2$ metric results in the lowest
residuals over short time horizons. The Stratonovich projection comes
a close second. Over medium time horizons, the It\^o-jet projection
out performs the It\^o-vector projection. The projection methods out-performed all
other methods like extended Kalman filter or assumed density filters. 

Second, in \cite{armstronglms} we compared the Hellinger residuals for different filters, where projection filters are w.r.t. the Hellinger metric. 
This second analysis indicates that the It\^o ADF and the It\^o-jet projection
are almost indistinguishable in their performance. A look at the
explicit formulae reveals that the difference between these two equations in the cubic sensor example with Gaussian densities  is of order $\epsilon^2$ whereas the difference between the other equations is of order only $\epsilon$.
Over the short term, the It\^o-vector projection gives the best
results. Over medium term, the It\^o-jet projection and the It\^o ADF give the best results.

We also note that in previous works such as \cite{brigo1,brigo2,armstrongbrigomcss} where we only studied the Stratonovich projection filter, filtering problems for systems like the cubic and quadratic sensors were studied. For such systems, the optimal filter density would often turn out to be bimodal and a projection filter based on a manifold consisting of mixtures of two Gaussians or of exponential families with fourth order polynomial exponents would track the optimal filter well, while approximated filters such as the Extended Kalman filter, Gaussian Assumed Density filters and even particle filters with the same number of parameters as the projection filters would fare worse than the projection filters in terms of $L^2$, Hellinger or L\'evy–Prokhorov norms of errors.  

\section{Conclusions and further work}
\label{section:conclusions}

The notion of projecting a vector field onto a manifold is unambiguous.
By contrast, there are multiple distinct generalizations of this notion to
SDEs, as summarized in \Cref{tab:projtypes}.

The two It\^o projections we introduced in this work can both be derived
from minimization arguments. However, the It\^o-jet projection has some clear advantages.
\begin{itemize}
\item  The It\^o-jet projection is the best approximation to the metric projection of the true solution and leads to a mean-squared error of order $O(t^2)$. By contrast, the It\^o-vector projection only tracks the true solution with an accuracy of $O(t)$ for the mean-square error.
\item  The It\^o-jet projection gives a more intuitive answer than
	   the It\^o-vector projection for the
	   low dimensional example of the cross-diffusion considered in \cite{armstronglms}.  
\item  The It\^o-jet projection gives better numerical results in the
	   longer term than the It\^o-vector projection in our application to filtering.	   
\item  The It\^o-jet projection has an elegant definition when written in terms
       of $2$-jets, which is described in \cite{armstronglms}. 
\end{itemize}

We have also seen that the Stratonovich projection satisfies an ad hoc minimization that is less appealing than the ones of the It\^o projections, since it requires a deterministic anchor point at time $0$. 
%
%
%
The It\^o-jet and  It\^o-vector projection arguments allow one
to derive new Gaussian approximations to non-linear filters, and new exponential and mixture filters more generally, although the more general cases have not been explored in \cite{armstronglms}. Some of the possibilities with different projections, metrics and manifolds are shown in Table \ref{tab:pfclass}. This could be investigated in further work to complete the table. 
In the Gaussian case we do explore in \cite{armstronglms} applying the methods summarized in this review, unlike previous Gaussian approximations to non-linear filters, the projection approximations
are derived by fully explicit minimization arguments rather than heuristic arguments.
Thus, the notion of projecting an SDE onto a manifold, coupled with information geometry, is able to give new results even for this well-worn topic of approximate nonlinear filtering.

A final investigation line could be in deriving approximations based on approximating bases that are not made of densities or their square roots. Working with densities has the advantage of allowing information geometry to act clearly, but at the same time puts strong constraints on the approximating bases. As a simple example, one might wish to use ``mixtures'' of Hermite polynomials, which are not densities, as a basis for the approximation. One might wish to investigate to what extent it is possible to use non-density bases while retaining an information geometric approach. 

\bigskip

\begin{table}[tbhp]
\begin{center}
\begin{tabular}{p{20mm}p{80mm}}\toprule
Projection  & Properties of drift term  \\ \midrule

It\^o-vector
&
(i)
Minimizes order 1 Taylor expansion (in $t$) of
norm of the expectation of the difference between $X$ \& $\phi(Y)$, namely $| E[X_t- \phi(Y_t)]|$.

(ii)
Given $B$ minimizing (but not zeroing) the $t$ term of the expansion for the mean square difference  $E[|X_t-\phi(Y_t)|^2]$, finds $A$ that minimizes the $t^2$ term while holding that $B$ fixed. 
\\ \addlinespace                     
 It\^o-jet
  &
Zeroes $t$ term and minimizes $t^2$ term of Taylor expansion of the mean square of the distance in  $\R^r$ or $M$  between $\pi(X)$ \& $\phi(Y)$, namely $\E [ d_M({\pi(X_t)},\phi(Y_t))^2 ]$ or  
$E [ |{\pi(X_t)}-\phi(Y_t)|_r^2 ]$.

\\ \addlinespace                     

Stratonovich
&
Similar to It\^o vector but for the Taylor series of the ``time-symmetric'' mean square difference between $X$ and its lower dimensional approximation $\phi(Y)$: 
\[\frac{1}{2} \left( \E [ | {X_{-t}} -  \phi(Y_{-t})|^2 ]  +
\E [ | {X_t} -  \phi(Y_t)|^2 ]  \right)\ \ \mbox{or} \ \  \]
\[ \frac{1}{2} \left( \E [ |{\pi(X_{-t})}-\phi(Y_{-t})|_r^2 ] +
\E [ |{\pi(X_t)}-\phi(Y_t)|_r^2 ] \right) \]
where negative time processes are defined ad hoc by propagating a second input Brownian motion backward in time. 
 \\
                      \bottomrule
\end{tabular}
\end{center}
\caption{ Projections and the associated optimality criteria}\label{tab:projtypes}
\end{table}

\newpage

\bibliographystyle{plain}
\bibliography{filtering}

\begin{thebibliography}{10}

\bibitem{aggarwal}
Nand~Lal Aggarwal.
\newblock Sur l'information de fisher.
\newblock In J.~Kamp{\'e}~de F{\'e}riet and C.~F. Picard, editors, {\em
  Th{\'e}ories de l'Information}, pages 111--117. Springer, Berlin Heidelberg,
  1974.

\bibitem{ahmedbook}
Nasir~Uddin Ahmed.
\newblock {\em Linear and Nonlinear Filtering for Scientists and Engineers}.
\newblock World Scientific, Singapore, 1998.

\bibitem{amari}
Shun'ichi Amari.
\newblock {\em Differential Geometrical Methods in Statistics}, volume~28 of
  {\em {L}ecture {N}otes in {S}tatistics}.
\newblock Springer Verlag, 1985.

\bibitem{armstrongBrigo}
John Armstrong and Damiano Brigo.
\newblock Stochastic filtering via {L}2 projection on mixture manifolds with
  computer algorithms and numerical examples.
\newblock {\em arXiv preprint arXiv:1303.6236}, 2013.

\bibitem{armstrongbrigoicms}
John Armstrong and Damiano Brigo.
\newblock Extrinsic projection of {I}t{\^o} {SDE}s on submanifolds with
  applications to non-linear filtering.
\newblock {\em In: Nielsen, F., Critchley, F., \& Dodson, K. (Eds),
  Computational Information Geometry for Image and Signal Processing, Springer
  Verlag}, 2016.

\bibitem{armstrongbrigomcss}
John Armstrong and Damiano Brigo.
\newblock {N}onlinear filtering via stochastic {PDE} projection on mixture
  manifolds in {$L^2$} direct metric.
\newblock {\em Mathematics of Control, Signals and Systems}, 28(1):1--33, 2016.

\bibitem{armstrongjetsgsi}
John Armstrong and Damiano Brigo.
\newblock Ito stochastic differential equations as 2-jets.
\newblock In {\em {G}eometric {S}cience of {I}nformation, {P}roceedings of the
  2017 {C}onference, {P}aris, {F}rance}, pages 543--551. Springer, 2017.

\bibitem{armstrongrspa}
John Armstrong and Damiano Brigo.
\newblock Intrinsic stochastic differential equations as jets.
\newblock {\em Proceedings of the Royal Society A: Mathematical, Physical and
  Engineering Sciences}, 474, 2018.

\bibitem{rossiferrucci}
John Armstrong, Damiano Brigo, and Emilio {Rossi Ferrucci}.
\newblock Projection of {SDE}s onto submanifolds.
\newblock {\em arXiv preprint arXiv:1810.03923}, 2018.

\bibitem{armstronglms}
John Armstrong, Damiano Brigo, and Emilio {Rossi Ferrucci}.
\newblock Optimal approximation of sdes on submanifolds: the ito-vector and
  ito-jet projections.
\newblock {\em Proceedings of the London Mathematical Society},
  119(1):176--213, 2019.

\bibitem{navigation}
Babak Azimi-Sadjadi and P.S. Krishnaprasad.
\newblock Approximate nonlinear filtering and its application in navigation.
\newblock {\em Automatica}, 41(6):945--956, 2005.

\bibitem{crisan}
Alan Bain and Dan Crisan.
\newblock {\em Fundamentals of stochastic filtering}, volume~3.
\newblock Springer, 2009.

\bibitem{barndorff}
Ole Barndorff-Nielsen.
\newblock {\em Information and Exponential Families}.
\newblock John Wiley \& Sons, New York, 1978.

\bibitem{belopolskaja}
Ya. Belopolskaja and Yu. Dalecky.
\newblock {\em Stochastic Equations and Differential Geometry}.
\newblock Mathematics and Its Applications, Vol. 30. Dordrecht, Kluwer Academic
  Publishers, Boston, London, 1990.

\bibitem{foi}
F.~Berefelt, Johan Hamberg, and John~W.C. Robinson.
\newblock Geometric aspects of nonlinear filtering.
\newblock {\em {T}ech. rep., {S}wedish {D}efence {R}esearch {A}gency}, 2003.

\bibitem{brigoscl1}
Damiano Brigo.
\newblock On the nice behaviour of the gaussian projection filter with small
  observation noise.
\newblock {\em Systems \& Control Letters}, 26(5):363--370, 1995.

\bibitem{brigoPHD}
Damiano Brigo.
\newblock {\em Filtering by Projection on the Manifold of Exponential
  Densities}.
\newblock Free University of Amsterdam, {P}h{D} dissertation, 1996.

\bibitem{brigoscl2}
Damiano Brigo.
\newblock New results on the gaussian projection filter with small observation
  noise.
\newblock {\em Systems \& Control Letters}, 28(5):273--279, 1996.

\bibitem{brigoIME}
Damiano Brigo and Bernard Hanzon.
\newblock On some filtering problems arising in mathematical finance.
\newblock {\em Insurance: Mathematics and Economics}, 22(1):53--64, 1998.

\bibitem{brigo1}
Damiano Brigo, Bernard Hanzon, and Fran{\c{c}}ois LeGland.
\newblock A differential geometric approach to nonlinear filtering: the
  projection filter.
\newblock {\em IEEE Transactions on Automatic Control}, 43(2):247--252, 1998.

\bibitem{brigo2}
Damiano Brigo, Bernard Hanzon, and Fran{\c{c}}ois LeGland.
\newblock Approximate nonlinear filtering by projection on exponential
  manifolds of densities.
\newblock {\em Bernoulli}, 5(3):495--534, 1999.

\bibitem{brigopistone}
Damiano Brigo and Giovanni Pistone.
\newblock Dimensionality reduction for measure-valued evolution equations in
  statistical manifolds.
\newblock In {\em Proceedings of the conference on Computational Information
  Geometry for Image and Signal Processing}. Springer, 2016.

\bibitem{brigopistone2}
Damiano Brigo and Giovanni Pistone.
\newblock Optimal approximations of the {F}okker-{P}lanck-{K}olmogorov
  equation: projection, maximum likelihood eigenfunctions and {G}alerkin
  methods.
\newblock {\em arXiv preprint arXiv:1603.04348}, 2016.

\bibitem{chaosfiltering}
Jochen Br\"ocker and Ulrich Parlitz.
\newblock Noise reduction and filtering of chaotic time series.
\newblock {\em {P}roc. {NOLTA} 2000}, 2000.

\bibitem{brzezniak}
Z.~Brze\'zniak and K.~D. Elworthy.
\newblock Stochastic differential equations on {B}anach manifolds.
\newblock {\em Methods Funct. Anal. Topology}, 6(1):43--84, 2000.

\bibitem{cipra}
Barry Cipra.
\newblock Engineers look to {K}alman {F}iltering for {G}uidance.
\newblock {\em {SIAM} {N}ews}, 26, 1993.

\bibitem{elworthy}
David Elworthy.
\newblock Geometric aspects of diffusions on manifolds.
\newblock In {\em {\'E}cole d'{\'E}t{\'e} de Probabilit{\'e}s de Saint-Flour
  XV--XVII, 1985--87}, pages 277--425. Springer, 1988.

\bibitem{geometryoffiltering}
David Elworthy, Yves~Le Jan, and Xue-Mei Li.
\newblock {\em The Geometry of Filtering}.
\newblock Frontiers in Mathematics. Birkhauser, 2010.

\bibitem{emery}
Michel Emery.
\newblock {\em Stochastic calculus in manifolds}.
\newblock Springer-Verlag, Heidelberg, 1989.

\bibitem{quantum2}
Qing Gao, Guofeng Zhang, and Ian~R. Petersen.
\newblock An exponential quantum projection filter for open quantum systems.
\newblock {\em Automatica}, 99:59--68, 2019.

\bibitem{gliklikh}
Yuri~E. Gliklikh.
\newblock {\em {Global and Stochastic Analysis with Applications to
  Mathematical Physics}}.
\newblock Theoretical and Mathematical Physics. Springer, London, 2011.

\bibitem{hanzonhut}
B.~Hanzon and R.~Hut.
\newblock New results on the {Projection} {Filter}.
\newblock In {\em {P}roceedings of the {E}uropean {C}ontrol {C}onference,
  {G}renoble, {F}rance}, pages 623--628. {ECC91}, 1991.

\bibitem{hanzon}
Bernard Hanzon.
\newblock A differential-geometric approach to approximate nonlinear filtering.
\newblock In C.T.J. Dodson, editor, {\em Geometrization of Statistical Theory},
  pages 219--223. ULMD Publications, University of Lancaster, 1987.

\bibitem{neuralencoding}
Yuval Harel, Ron Meir, and Manfred Opper.
\newblock A tractable approximation to optimal point process filtering:
  Application to neural encoding.
\newblock In C.~Cortes, N.~Lawrence, D.~Lee, M.~Sugiyama, and R.~Garnett,
  editors, {\em Advances in Neural Information Processing Systems}, volume~28.
  Curran Associates, Inc., 2015.

\bibitem{hazewinkel}
Michiel Hazewinkel, Steven~I. Marcus, and H.J. Sussmann.
\newblock Nonexistence of finite-dimensional filters for conditional statistics
  of the cubic sensor problem.
\newblock {\em Systems \& control letters}, 3(6):331--340, 1983.

\bibitem{hsu}
Elton~P. Hsu.
\newblock {\em Stochastic Analysis on Manifolds}.
\newblock Contemporary Mathematics. American Mathematical Society, 2002.

\bibitem{jazwinski}
Andrew~H. Jazwinski.
\newblock {\em Stochastic Processes and Filtering Theory}.
\newblock Academic Press, New York, 1970.

\bibitem{soatto}
Eagle~S. Jones and Stefano Soatto.
\newblock Visual-inertial navigation, mapping and localization: A scalable
  real-time causal approach.
\newblock {\em The International Journal of Robotics Research}, 30(4):407--430,
  2011.

\bibitem{kushner}
Harold~J. Kushner.
\newblock Approximations to optimal nonlinear filters.
\newblock {\em Automatic Control, IEEE Transactions on}, 12(5):546--556, 1967.

\bibitem{circular}
Anna Kutschireiter, Luke Rast, and Jan Drugowitsch.
\newblock Projection filtering with observed state increments with applications
  in continuous-time circular filtering.
\newblock {\em {IEEE} Transactions on {S}ignal {P}rocessing}, 70, 2022.

\bibitem{ocean}
Pierre~F.J. Lermusiaux.
\newblock Uncertainty estimation and prediction for interdisciplinary ocean
  dynamics.
\newblock {\em Journal of Computational Physics}, 217(1):176--199, 2006.

\bibitem{liptser}
R.S Liptser and A.N. Shiryayev.
\newblock {\em Statistics of Random Processes {I}, General Theory}.
\newblock Springer Verlag, Berlin, 1978.

\bibitem{hazardposition}
Yan Ma, Yu-xin Zhao, Li-juan Chen, and Chang Shuai.
\newblock Hazard position estimation of sudden accidents based on projection
  filter.
\newblock {\em Systems Engineering - Theory \& Practice}, 35(3):651, 2015.

\bibitem{maybeck}
Peter~S. Maybeck.
\newblock {\em Stochastic models, estimation, and control}, volume~3.
\newblock Academic press, 1982.

\bibitem{newton1}
Nigel~J. Newton.
\newblock An infinite-dimensional statistical manifold modelled on {H}ilbert
  space.
\newblock {\em {J}. {F}unct. {A}nal.}, 263(6):1661--1681, 2012.

\bibitem{newton2}
Nigel~J. Newton.
\newblock Information geometric nonlinear filtering.
\newblock {\em {I}nfin. {D}imens. {A}nal. {Q}uantum {P}robab. {R}elat. {T}op.},
  2(18):1550014, 24, 2015.

\bibitem{pardoux}
{\'E}tienne Pardoux.
\newblock Filtrage non lin{\'e}aire et {\'e}quations aux d{\'e}riv{\'e}es
  partielles stochastiques associ{\'e}es.
\newblock In {\em Ecole d'Et{\'e} de Probabilit{\'e}s de Saint-Flour XIX 1989},
  pages 68--163. Springer, 1991.

\bibitem{pistoneannals}
Giovanni Pistone and Carlo Sempi.
\newblock An infinite-dimensional geometric structure on the space of all the
  probability measures equivalent to a given one.
\newblock {\em Ann. Statist.}, 23(5):1543--1561, 10 1995.

\bibitem{rao45}
C.~R. Rao.
\newblock Information and accuracy attainable in the estimation of statistical
  parameters.
\newblock {\em {B}ull. {C}alcutta {M}ath. {S}oc.}, pages 81--91, 1945.

\bibitem{rogerswilliams}
L.~C.~G. Rogers and David Williams.
\newblock {\em Diffusions, {M}arkov processes, and martingales. {V}ol. 2: {I}to
  calculus}.
\newblock Cambridge University Press, 1987.

\bibitem{MLEPartial}
Simone Surace and Jean-Pascal Pfister.
\newblock Online maximum likelihood estimation of the parameters of partially
  observed diffusion processes.
\newblock {\em IEEE Transactions on Automatic Control}, 64, 09 2017.

\bibitem{arbitrarylike}
Filip Tronarp and Simo S\"arkk\"a.
\newblock Updates in bayesian filtering by continuous projections on a manifold
  of densities.
\newblock In {\em ICASSP 2019 - 2019 IEEE International Conference on
  Acoustics, Speech and Signal Processing (ICASSP)}, pages 5032--5036, 2019.

\bibitem{quantum1}
Ramon van Handel and Hideo Mabuchi.
\newblock Quantum projection filter for a highly nonlinear model in cavity
  {QED}.
\newblock {\em Journal of {O}ptics {B}: {Q}uantum and {S}emiclassical
  {O}ptics}, 7, 2005.

\bibitem{changepoint}
M.~H. Vellekoop and J.~M.~C. Clark.
\newblock A nonlinear filtering approach to changepoint detection problems:
  Direct and differential-geometric methods.
\newblock {\em SIAM Review}, 48(2):329--356, 2006.

\bibitem{fiber}
Xian~Miao Zhang, Rong~Wu Wang, Hai~Bo Wu, and Bugao Xu.
\newblock Automated measurements of fiber diameters in melt-blown nonwovens.
\newblock {\em Journal of Industrial Textiles}, 43(4):593--605, 2014.

\end{thebibliography}

\end{document}